\documentclass[runningheads]{llncs}

\usepackage{amssymb}
\usepackage{graphicx}
\usepackage{multirow}
\usepackage{float}
\usepackage{amsmath}
\usepackage[T1]{fontenc}
\usepackage{url}

\newcommand{\BG}{\emph{BG}}

\newcommand{\transpose}[1]{{#1}^{t}}

\newcommand{\rank}{\mathrm{\mathop{rank}}}

\newcounter{algorithm}
\newcommand{\algorithmName}[1]{%
\vspace{10pt}
\hrule height 1pt
\vspace{3pt}
\refstepcounter{algorithm}
\noindent {\bf Algorithm \arabic{algorithm}:} #1
\vspace{3pt}
\hrule
\vspace{2pt}}
\newcommand{\algorithmInput}[1]{\noindent{\bf Input:} #1 \\}
\newcommand{\algorithmOutput}[1]{\noindent{\bf Output:} #1 \\}
\newcommand{\algorithmComplexity}[1]{\noindent{\bf Complexity:} #1 \\}
\newcommand{\algorithmBegin}{\noindent{\bf Procedure:} \begin{enumerate}}
\newcommand{\algorithmEnd}{\end{enumerate}}
\newcommand{\algorithmRemark}[1]{\noindent{\bf Remark:} #1 \\}
\newcommand{\algorithmLine}{%
\hrule height 1pt
\vspace{10pt}
}

\newcommand{\True}{\emph{True}}
\newcommand{\False}{\emph{False}}

\def\BB{$\cal B$}

\newcommand\FALSE{\emph{False\/}}

\newcommand\ldotspec{\ldots\,}

\def\yy#1{\hbox to 4mm{\hfil#1}}
\def\zz#1{\hbox to 6mm{\hfil#1}}
\def\zzz#1{\hbox to 8.5mm{\hfil#1}}

\def\tablewithrules#1\par{
\noindent
  \vbox{
     \hbox to\hsize{
        \hfil
        \vbox{
          \halign {
              ##\hfil&&\quad##\hfil\cr #1
          }
        }
        \hfil
     }
  }
}
\begin{document}

\mainmatter  % start of an individual contribution

\title{Implementation of a Unimodularity Test}

\titlerunning{Unimodularity Test}

\author{Matthias Walter$^1$
\and Klaus Truemper$^2$}

\authorrunning{M. Walter and K. Truemper}

\institute{$^1$Institute of Mathematical Optimization,\\
           University of Magdeburg
              {\lq\lq}Otto von Guericke{\rq\rq},\\
           39106 Magdeburg, Germany \\
           $^2$Department of Computer Science,\\
           University of Texas at Dallas,\\
           Richardson, TX 75083, U.S.A.}

\maketitle

%\noindent
%Corresponding author:\\
%M.~Walter \\
%walter@www.math.uni-magdeburg.de

\medskip

\begin{abstract}

This paper describes implementation and computational
results of a polynomial test of total unimodularity. 
The test is a simplified
version of a prior method. The program also
decides two related unimodularity properties.
The software is available free of charge
in source code form under the Boost Software License.

\end{abstract}


%-----------------------------------------------
\section{Introduction}
\label{sec:Introduction}
%-----------------------------------------------

This paper describes the implementation of a simplified version
of the polynomial test \cite{Truemper90} for total unimodularity.
The program also decides absence/presence of
two related types of unimodularity.
The computer program is available free of charge in source code
from two sites \cite{Software}
under the Boost Software License \cite{Boostlicense}.
Computational effectiveness is demonstrated for nontrivial
test instances.

We begin with a well-known definition.
An integer matrix $A$ is \emph{totally unimodular} (t.u.)
if every square submatrix $D$
of $A$ has $\det D = 0$ or $\pm1$.
This property was introduced by
\cite{Hoffman56} with slightly different terminology.
The reference establishes a key result for t.u. matrices:
The inequality $Ax \le b$ has all basic solutions integer for
all integer vectors $b$ if and only if $A$ is t.u.

Several other concepts are closely related to total unimodularity;
see, for example,
\cite{Hoffman78,%
Truemper78,%
TruemperRC78,%
Truemper80%
}.
Here, we cover the following two properties.
An integer matrix $A$
is \emph{unimodular} if for every column basis $C$
of $A$, the maximal square submatrices $C^i$ of $C$
satisfy $\gcd_i\det C^i = 1$.
In the special case of unimodularity
where each $C^i$ has $\det C^i = 0$ or $\pm1$, the
matrix $A$ is \emph{strongly unimodular}.
Analogously to the key result for t.u. matrices,
the equation  $Ax = b$ has all basic solutions integer for
all integer vectors $b$ if and only if $A$ is unimodular.
The next theorem is taken from
\cite{Truemper78}.

\newpage

\begin{theorem}
\label{thm:unimodularity}
Let $A$ be an integer matrix.
\begin{itemize}
\item[\rm 1.] $A$ is unimodular if and only if,
for an arbitrarily selected
column basis matrix $C$ of $A$, $\gcd_i\det C^i = 1$
and the solution $X$ of $CX = A$ is t.u.

\item[\rm 2.] $A$ is strongly unimodular if and only if
both $A$ and its transpose $A^t$ are unimodular.
\end{itemize}
\end{theorem}

The question whether
a given column basis $C$ of an integer matrix $A$
satisfies $\gcd_i\det C^i = 1$ can be efficiently answered
via the Smith Normal Form of \cite{Smith1861}.
Thus, testing for any of the above properties is readily
reduced to testing for total unimodularity. The software
carries out the same reductions.
We skip implementation details and
focus on the test of total unimodularity.

Let $A$ be a given integer matrix and $B$ be the binary matrix
derived from $A$ by converting each $-1$ to $1$. Define $I$
to be the identity of appropriate order. It is well known that
the following four steps
decide absence/presence of total unimodularity.

In the first step, a trivial check verifies
that all nonzeros of $A$ are
$1$ or $-1$. Clearly, $A$ is not t.u.
if the matrix fails this test. 
In the second
step, it is checked whether the binary matroid $M(B)$ represented
by the matrix $[I|B]$ has the property of
\emph{regularity} defined by \cite{Tutte58}. For our purposes,
it is convenient to declare $M(B)$
to be regular if $B$ can be signed
to become t.u. Testing for regularity of $M(B)$ is the most
difficult step. If the answer is
negative, then $A$ cannot be t.u. So assume that $M(B)$
is regular.
In the third step,
suitable signing of the 1s of $B$ converts that matrix
to a t.u. matrix $A'$. This step
is based on the uniqueness of such signing, up to scaling,
proved by \cite{Camion63}. The process is quite straightforward.
The uniqueness result of \cite{Camion63}
is once more employed in the fourth step,
where it is checked whether $A'$ can by column and row scaling
be converted to $A$. This test is very easy.
The matrix $A$ is t.u. if and only if such scaling is possible.

As far as we know, every prior polynomial algorithm for
testing matroid regularity uses the regular matroid decomposition of
\cite{Seymour80}. Indeed, implicit in the cited reference
is already one such scheme, provided the proofs are implemented
in suitable algorithmic steps.
Of the prior methods \cite{Bixby86,Cunningham78,Truemper90},
the scheme of \cite{Truemper90} has lowest order, which for
a binary $m\times n$ matrix is $O((m+n)^3)$. We use the latter
method as basis for the implementation. We say ``as basis'' since,
for a first attempt, full implementation of all features of
\cite{Truemper90} was rather daunting.
Hence, we opted for a simplified version
that avoids complex operations but still is in
the spirit of the method. The next section summarizes
that implementation while pointing out differences
to \cite{Truemper90}. Section~\ref{sec:EnumerativeMethods}
introduces two na\"{\i}ve enumerative tests 
for comparison purposes. 
Section~\ref{sec:ComputationalResults}
describes computational results for nontrivial test instances.
Section~\ref{sec:DetailsofImplementation}
contains technical details of the implementation.

%-----------------------------------------------
\section{Summary of Implementation}
\label{sec:SummaryofImplementation}
%-----------------------------------------------

This section describes the implementation and compares it with
the method of \cite{Truemper90}. For an abbreviated terminology,
\emph{implemented method} refers to the method implemented to-date,
while \emph{original method} is the scheme of the cited reference.

We strive for an intuitive discussion so that
salient ideas emerge and a clutter of technical details is avoided.
The reader not familiar with matroid theory may want to rely
on the introductory discussion of binary,
graphic, and regular matroids of \cite{Truemper98}
before proceeding.

We first review key concepts for testing
matroid regularity. For details, see \cite{Seymour80,Truemper98}.

%------------------------------------
\subsection{Key Concepts}
\label{sec:KeyConcepts}
%------------------------------------

We begin with some definitions regarding matrix notation.
Let $A$ be a matrix whose rows (resp.\ columns) are indexed
by a set $Z$ (resp.\ $Y$). For any subsets
$Z'\subseteq Z$ and $Y'\subseteq Y$, the submatrix
of $A$ indexed by $Z'$ and $Y'$ is denoted by
$A_{Z',Y'}$.
If $Z'$ is a singleton set, say consisting of an
element $z$, then we use just $z$ instead of
$\{z\}$ in the above notation. A singleton set containing
an element $y$ of $Y$ is handled analogously. In particular,
$A_{z,y}$ is the entry of $A$ indexed by
$z$ and $y$. The abbreviated notation is not ambiguous
since we always employ upper-case letters for sets and lower-case
letters for elements of sets.

For a binary matrix $[I|B]$,
let $X$ and $Y$ index the columns of the submatrices $I$ and $B$,
respectively. The matroid
$M(B)$ has $X\cup Y$ as groundset.
For $X'\subseteq X$ and $Y'\subseteq Y$,
the subset $X'\cup Y'$ of $X\cup Y$ is independent
in the matroid if the column submatrix of
$[I|B]$ indexed by $X'\cup Y'$ has linearly independent columns.
We avoid explicit display of the
matrix $I$ by indexing the rows of
$B$ by $X$, in addition to the column index set $Y$.
Then $X'\cup Y'$ is independent in $M(B) $ if and only if the
submatrix $B'= B_{X \setminus X',Y'}$ has independent columns.

We consider two matrices \emph{equal} if they become
numerically the same
under suitable row and column permutations.
The indices of rows and columns are ignored in the comparison.
It is convenient that we apply matroid terminology for $M(B)$
to $B$ as well. Thus, $B$ is \emph{regular} if $M(B)$ has that
property, that is, if $B$ can be signed to become a t.u. matrix.

For $x\in X$ and $y\in Y$,
a \emph{pivot} on a nonzero entry $B_{x,y}$
of the matrix $[I|B]$ is the
customary set of elementary row operations.
In the reduced notation where $I$ is not explicitly listed,
the pivot converts $B$ to a matrix $B'$ that agrees
numerically with $B$ except for the entries $B'_{i,j}$
where $i\neq x$, $j\neq y$, and
$B_{x,j} = B_{i,y} = 1$. The matrix $B'$ has the same index sets
as $B$ except that the indices $x$ and $y$ have traded places.
The matrix $B$ is regular if and only if this holds for $B'$.

$M(B)$ and $B$ are \emph{graphic}
if there exists an undirected graph $G$
with edges indexed by the elements of $X\cup Y$ such that
the edge sets of subgraphs of $G$ without any cycle are precisely
the independent sets of $M(B)$. Note that each zero column of
$B$ corresponds to a loop of $G$.
$M(B)$ and $B$ are \emph{cographic}
if the transpose of $B$, denoted by $B^t$, is graphic.
$M(B)$ and $B$ are \emph{planar} if $B$ is graphic and cographic.

A graphic $B$ is regular.
Since $B$ can be signed to become t.u. if and only if
this is so for $B^t$, a cographic $B$ is regular as well.
There exist very fast algorithms for deciding
whether $B$ is graphic \cite{Bixby80,Fujishige80}.

Define $\BG(B)$ to be the bipartite graph with node set $X\cup Y$
where an undirected edge joins nodes $x\in X$ and $y\in Y$
if $B_{x,y} = 1$. Declare the matrix $B$ to be \emph{connected} if
the graph $\BG(B)$ is connected.
Define the \emph{length} of $B$, denoted by $s(B)$, to be the
number of rows plus the number of columns of $B$.
We allow matrices to have no rows or columns. The rank
of any such matrix is 0.

If $B$ has zero or unit vector rows or columns,
or has duplicate rows or
columns, then recursive deletion of zero/unit vectors
and of duplicates
except for representatives,
reduces $B$ to a \emph{simple} matrix $B'$.
The matrix $B$ is regular
if and only if this is so for $B'$.
If $B'$ has no rows or columns, then $B$ is regular.
For the definitions to follow, we assume that $B$ is simple.

If matrix $B$ has block structure, say with blocks $B^i$,
$i=1$, $2,\ldotspec$, then
$B$ is a \emph{$1$-sum} of the blocks $B^i$,
and $B$ is regular if and only if each of the
blocks $B^i$ has that property.
The implemented method detects blocks via $\BG(B)$.
For the remaining discussion of this section, we assume
that the simple $B$ has no such block structure,
which is equivalent to assuming that $B$ is connected.

Suppose $B$ has the form

\begin{figure}[H]
\begin{center}
\includegraphics[width=40mm]{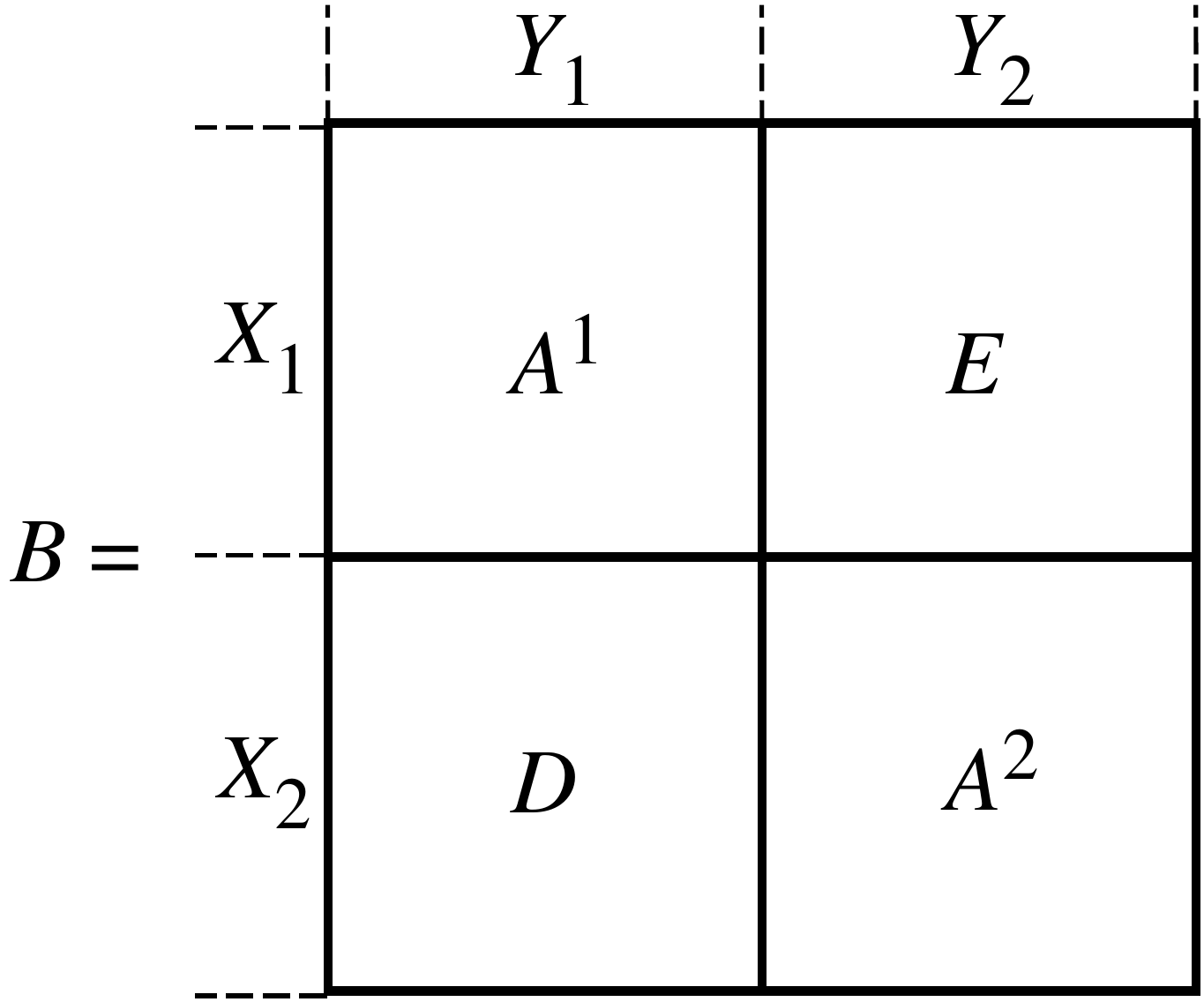}%
\caption{Separation of $B$}%
\label{fig:klseparation}
\end{center}
\end{figure}

\noindent
Define $k=\rank(D)+\rank(E)+1$. If the lengths of $A^1$ and $A^2$
satisfy $s(A^1)\ge k$ and $s(A^2)\ge k$, then $B$ has
a \emph{$k$-separation}. If either
$s(A^1)= k-1$ and $s(A^2)\ge k$, or
$s(A^1)\ge k$ and $s(A^2)=k-1$, then $B$ has a
\emph{deficient $k$-separation}. If, for some $l\ge k$,
$s(A^1)\ge l$ and $s(A^2)\ge l$, then $B$ has
a \emph{$(k|l)$-separation}.

If the submatrix $E$ of $B$ is nonzero, then by pivots in $E$
we can always obtain a matrix $B'$ with the same type of
separation where $E' = 0$ and $\rank(D') = k-1$. In the discussion
below, we assume that $B$ itself is of that form.

In the case of a 2-separation, the submatrix $D$
of Fig.~\ref{fig:klseparation} has rank equal to 1,
and $B$ has the following form.

\begin{figure}[H]
\begin{center}
\includegraphics[width=40mm]{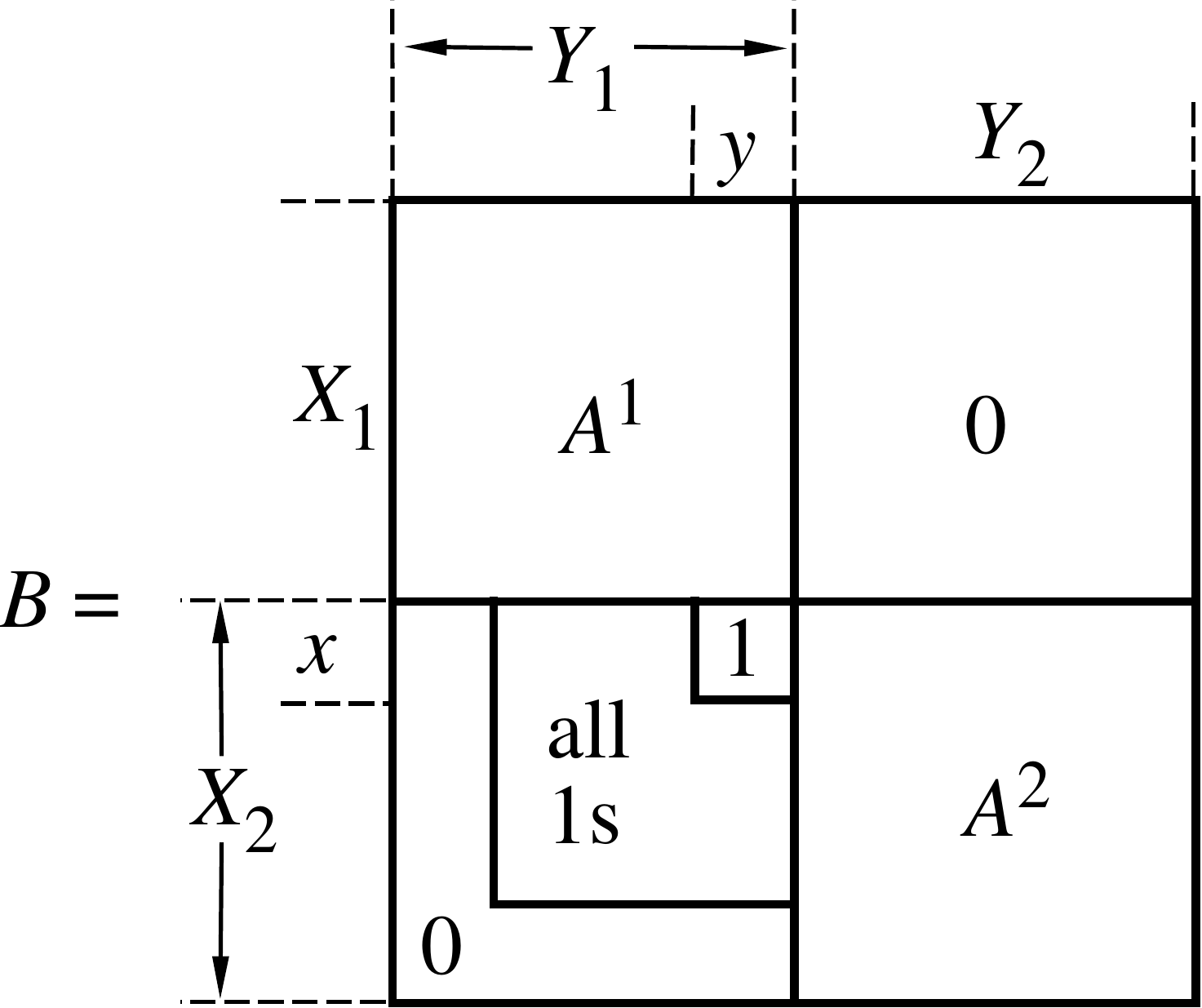}%
\caption{$2$-Sum Case}%
\label{fig:2sumCase}
\end{center}
\end{figure}

\noindent
If in addition both $A^1$ and $A^2$ have at least one entry,
then $B$ is a \emph{$2$-sum} with the following
\emph{component matrices} $B^1$ and $B^2$.

\begin{figure}[H]
\begin{center}
\includegraphics[width=63mm]{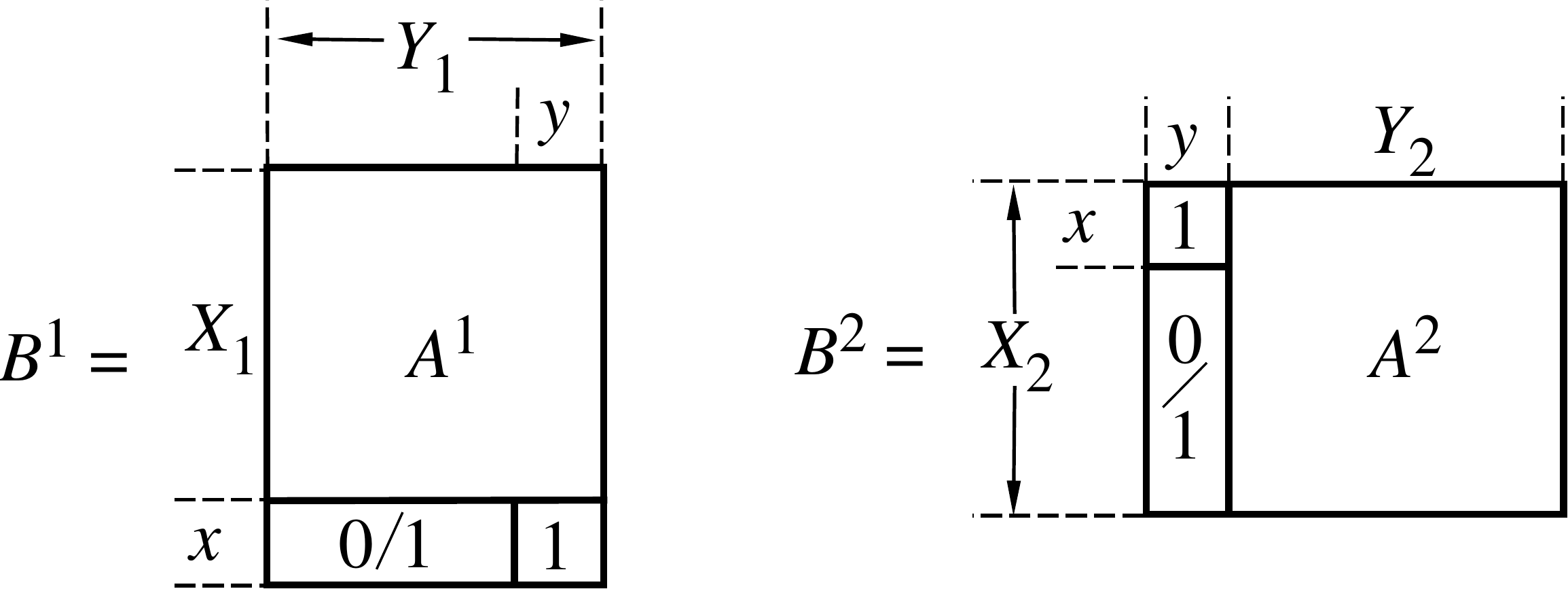}%
\caption{$2$-sum Components}%
\label{fig:2sumComponents}
\end{center}
\end{figure}

\noindent
Let $\overline{B}$ be the submatrix $B_{x,y}$ of $B$ in
Fig.~\ref{fig:2sumCase}. The same submatrix occurs in
$B^1$ and $B^2$ of Fig.~\ref{fig:2sumComponents}.
That submatrix has rank equal
to 1. It is called
the \emph{connecting} submatrix of the 2-sum decomposition.
When $B^1$ and $B^2$ are overlaid such that the two
connecting submatrices are identified, then straightforward
computations produce the matrix $B$ of Fig.~\ref{fig:2sumCase};
see Chapter~8 of \cite{Truemper98}.
The matrix $B$ of Fig.~\ref{fig:2sumCase}
is regular if and only if this is so for $B^1$ and
$B^2$.

Continuing the discussion of matrix features,
suppose that $B$ is connected, simple, and has no 2-sum
decomposition. Such a matrix is called \emph{$3$-connected}.
Assume that a 3-connected $B$ has a $(3|l)$-separation
where $l\ge 4$.
It is not difficult to prove that pivots can produce in $B$
the following structure for
the submatrices $A^1$, $A^2$, and $D$, where
$\overline{D}$ has the same rank as $D$,
that is, 2.

\begin{figure}[H]
\begin{center}
\includegraphics[width=40mm]{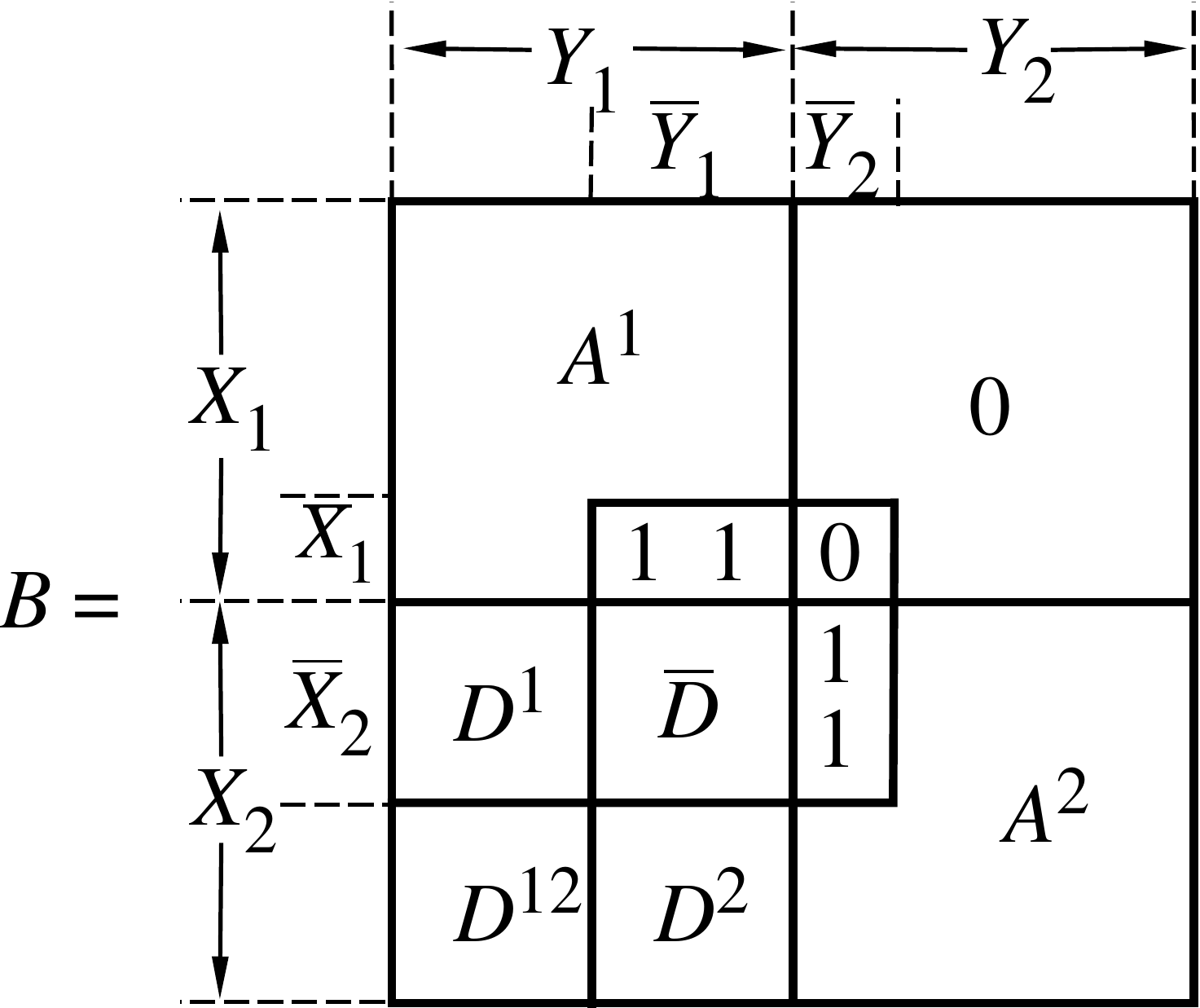}%
\caption{$3$-sum Case}%
\label{fig:3sumCase}
\end{center}
\end{figure}

\noindent
Then $B$ is a \emph{$3$-sum} with the following
\emph{component matrices}
$B^1$ and $B^2$.

\begin{figure}[H]
\begin{center}
\includegraphics[width=74mm]{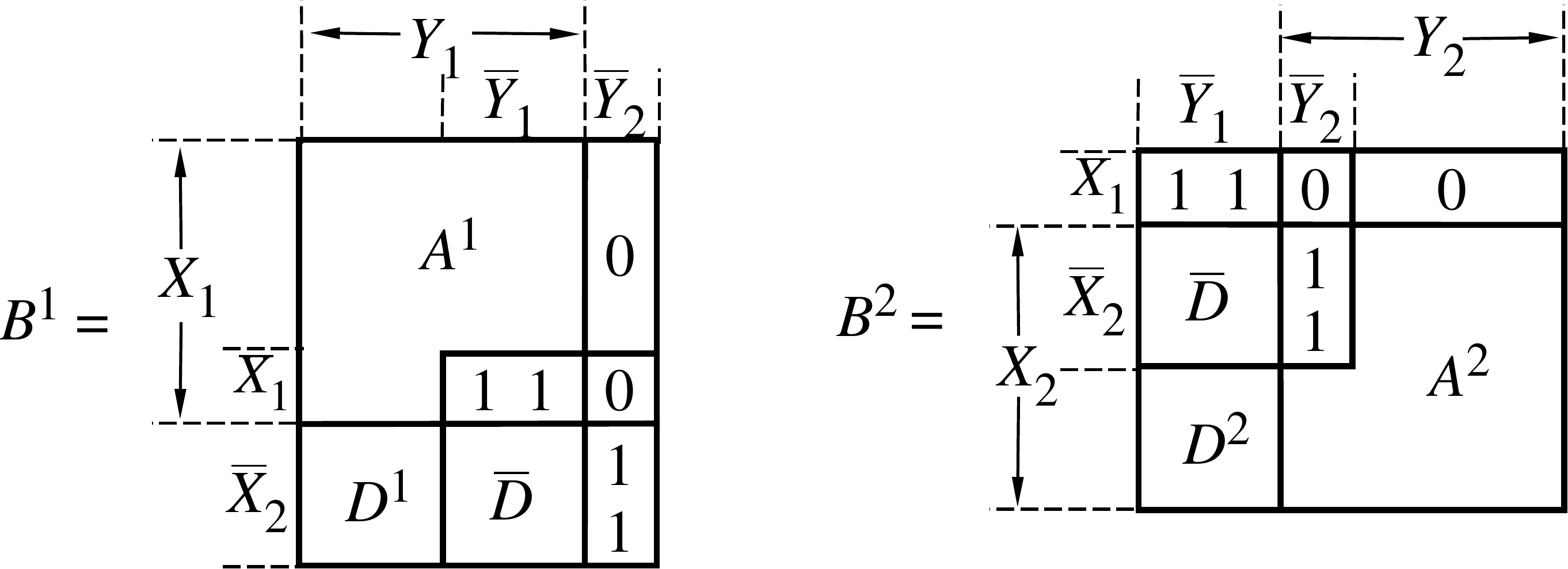}%
\caption{$3$-sum Components}%
\label{fig:3sumComponents}
\end{center}
\end{figure}

\noindent
For $B$ of Fig.~\ref{fig:3sumCase},
define $\overline B$ to be the submatrix
$B_{\overline X_1 \cup \overline X_2,
     \overline Y_1 \cup \overline Y_2}$.
The same submatrix occurs in
$B^1$ and $B^2$ of Fig.~\ref{fig:3sumComponents}.
Evidently, $\overline B$ is a $3\times 3$ matrix, so
$\overline{D}$ is a $2\times 2$ matrix. Since
$\overline{D}$ has rank equal to 2, it must be an identity
matrix or have exactly three 1s.
The matrix $\overline B$ is the \emph{connecting} submatrix
of the $3$-sum decomposition. It is easy to check
that $\overline B$ is graphic. The corresponding graph
is the wheel $W_3$ with three spokes.
When $B^1$ and $B^2$ are overlaid such that the two
connecting submatrices are identified, then, analogously
to the 2-sum case, straightforward
computations produce the matrix $B$ of Fig.~\ref{fig:3sumCase}.

The matrix $B$ of Fig.~\ref{fig:3sumCase}
is regular if and only if this is so for $B^1$ and
$B^2$ of Fig.~\ref{fig:3sumComponents}.
The above 2-sums and 3-sums decompositions can be
found by the matroid intersection algorithm of
\cite{Edmonds65} plus some pivots.

Finally, there is a regular matroid on 10 elements called
$R_{10}$. There are only two matrices that represent $R_{10}$.
They are $B^{10.1}$ and $B^{10.2}$ below.

\begin{figure}[H]
\begin{center}
\includegraphics[width=74mm]{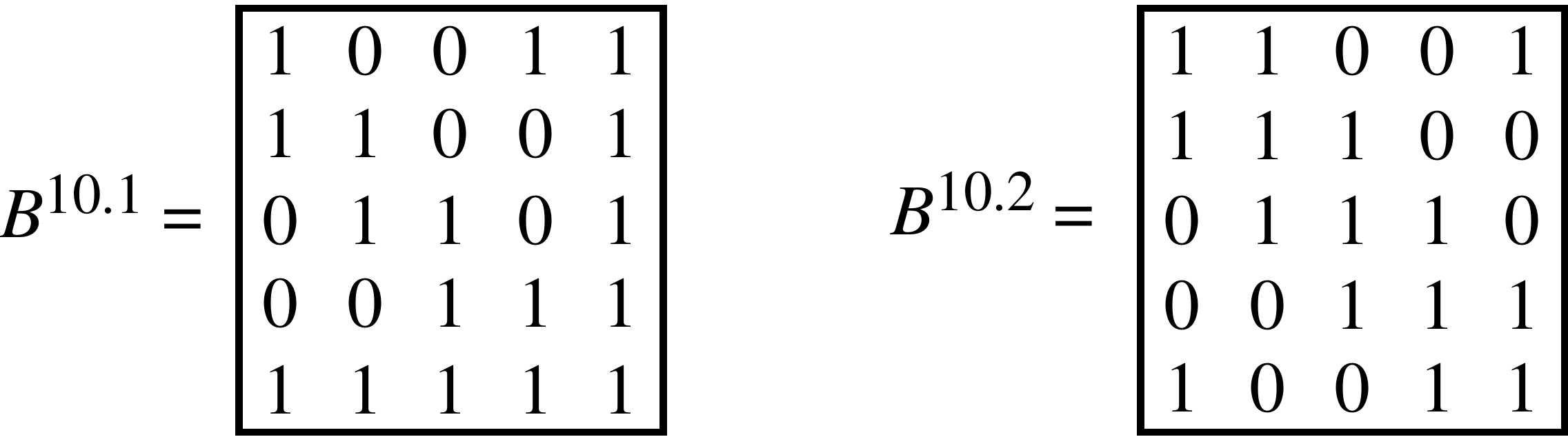}%
\caption{Matrices $B^{10.1}$ and $B^{10.2}$ for $R_{10}$}%
\label{fig:MatricesforR10}
\end{center}
\end{figure}

\noindent
$R_{10}$ is the smallest regular matroid that is
not graphic and not cographic. It is simple, connected, and does not
have a 2- or 3-sum decomposition.

In slightly different form and distributed among a number
of references,
the concepts and ideas stated above were known prior
to 1978. But that knowledge was not sufficient
to establish a polynomial testing algorithm for regularity.
That situation changed in 1978
when Seymour constructed the decomposition theorem for the regular
matroids \cite{Seymour80}, which supports efficient testing
of regularity.
Indeed, the cited reference implicitly already contains
such a scheme, provided certain nonconstructive proofs are
replaced by constructive ones involving polynomial subroutines.
A simplified version of the theorem that suffices
for present purposes is stated next.

\begin{theorem}
\label{thm:Seymour}
For any regular matrix, at least one of the statements
{\rm (i)}-{\rm (vi)} applies.
\begin{itemize}
\item[\rm (i)]   $B$ is graphic or cographic.
\item[\rm (ii)]  $B$ has a zero or unit vector row or column,
                     or has duplicate rows or columns.
\item[\rm (iii)] $B$ is simple and a $1$-sum.
\item[\rm (iv)]  $B$ is simple, connected, and a $2$-sum.
\item[\rm (v)]   $B$ is 3-connected and, after suitable pivots,
                     has the form of a $3$-sum that corresponds
                     to a $(3|6)$-separation. Let
                     $\overline B$ with row index
                     set $\overline X$ and column index set 
                     $\overline Y$ be any
                     $3$-connected nongraphic and noncographic
                     submatrix of $B$. Then there is a
                     $(3|6)$-separation of $B$, say defined by
                     index sets $X_1$, $X_2$, $Y_1$, and $Y_2$,
                     such that 
                     $X_1 \cap \overline X$,
                     $X_2 \cap \overline X$,
                     $Y_1 \cap \overline Y$, and
                     $Y_2 \cap \overline Y$
                     define a $(3|6)$-separation of $\overline B$.
\item[\rm (vi)]  $B$ is equal to $B^{10.1}$ or $B^{10.2}$ of
                     {\rm Fig.~\ref{fig:MatricesforR10}}.
\end{itemize}
\end{theorem}

Theorem~\ref{thm:Seymour} supports the following polynomial algorithm
for testing regularity. Given a matrix $B$, check with
one of the methods of \cite{Bixby80,Fujishige80} whether $B$ or
$B^t$ is graphic. If this is so, $B$ has been proved to be regular.
Otherwise, reduce $B$ to a simple matrix and check if it has a
1-, 2-, or 3-sum decomposition, in that order, using the
graph $\BG(B)$ for the 1-sum case and
the matroid intersection algorithm of \cite{Edmonds65}
for the 2- and 3- sum cases.
If a decomposition is detected, carry out the decomposition
and apply the algorithm recursively to the components.
Otherwise, check if
$B$ is equal to one of the matrices $B^{10.1}$ and $B^{10.1}$
of Fig.~\ref{fig:MatricesforR10}. If this is the case,
$B$ is regular. Otherwise, declare $B$ to be nonregular.

References \cite{Bixby86,Cunningham78} describe sophisticated
versions of the above approach, with bound
$O((m+n)^{4.5})(\log(m+n))^{0.5})$ for \cite{Bixby86}
and $O((m+n)^5)$ for \cite{Cunningham78}.
A lower complexity can be achieved when
the tests of graphicness and the search for decompositions
are intertwined, and when the latter search is carried out
by a certain induced decomposition scheme instead of the
matroid intersection algorithm \cite{Edmonds65}.
This is done in \cite{Truemper90}, producing
a test with bound $O((m+n)^3)$. The next section gives 
insight into that algorithm, which according to the convention
introduced earlier is called the \emph{original method}.
In the description below, the emphasis is on providing intuitive
insight into the method instead of a mathematically precise
specification, which is included in \cite{Truemper90}.

%------------------------------------
\subsection{Original Method}
\label{sec:OriginalMethod}
%------------------------------------

The method initializes a set \BB\ with a matrix $B_0$ that is to be
tested for regularity. The method removes
(resp.\ adds) matrices from (resp.\ to)
\BB\ until \BB\ becomes empty or nonregularity
of some matrix in \BB\ has been proved.
In the former (resp.\ latter) case,
$B$ has been proved to be regular (resp.\ not regular).
For some matrices of \BB, additional information is recorded.
Details are covered in the description of the method.

\begin{itemize}
\item[0.] If \BB\ is empty, declare the initial matrix $B_0$
to be regular, and stop.
Otherwise remove an arbitrary matrix $B$ from \BB.
If $B$ is known to be 3-connected, go to Step~3.

\item[1.]
If $B$ is not simple, remove zero and unit vectors
and reduce duplicate vectors to
representatives. If $B$ is not connected,
carry out a 1-sum decomposition via
the graph $\BG(B)$, place the components into \BB,
and go to Step~0.

\item[2.]
Determine a sequence of
\emph{nested 3-connected submatrices}
$N^1,\ldotspec$, $N^k $ for some $k\ge 1$, where
(1) $N^1$ is graphic and the corresponding graph $G$ is 
the wheel $W_3$,
(2) for each $i>1$, $N^i$ contains $N^{i-1}$ as proper submatrix
and the lengths of $N^i$ and $N^{i-1}$
satisfy $s(N^i) \le s(N^{i-1}) + 3$. We skip detailed discussion of
the steps finding such a sequence. They are described in
procedure FIND-B and EXTEND-B of the original method,
except that FIND-B or EXTEND-B assume 3-connectedness of $B$
and can be trivially modified to detect 2-sum decompositions.
Suffice it to say here that the steps make repeated use of
breadth-first-search (BFS) in the graph $\BG(B)$ or a closely related
graph, and of certain \emph{path shortening pivots}.
If the modified FIND-B or EXTEND-B detects a 2-sum decomposition,
 place the components of the 2-sum into
\BB\ and return to Step~0. Otherwise,
the last matrix $N^k$ of the sequence of
nested 3-connected submatrices is equal to $B$.
Place $B$ into \BB, suitably
record with it the 3-connected extension sequence,
and go to Step~0.

\item[3.] If $B$ is not supplied with a sequence of
nested 3-connected submatrices, construct such a sequence
as described in Step~2,
except that the 2-sum decomposition case cannot occur.
Regardless of the case,
each 3-connected submatrix $N^i$ triggers
additional testing as described next.
Recall that $N^1$ is graphic, indeed
represents the wheel $W_3$. It is assumed
inductively that $N^{i-1}$ is graphic and possibly planar.

\item[4.] If $N^{i-1}$ is planar, check whether $N^i$ is planar.
If the answer is negative, check
if $N^i$ is graphic or cographic; if $B$ turns out
to be cographic,
apply the transpose operator to $B$ and its submatrices.
 Thus, there are three possible outcomes:
$N^i$ is planar, or graphic but not
cographic, or not graphic and not cographic.
The test is carried out by TEST-C of the original method.
It is very efficient due a key result of \cite{Whitney33}
according to which a 3-connected graphic matrix
has exactly one corresponding graph.
If $N^i = B$ and $N^i$ is planar or graphic, then
$B$ is regular; go to Step~0.

\item[5.] Determine
whether any one of certain $(3|l)$-separations, $l\ge 3$, of
$N^i$ can be extended to
 a $(3|l')$-separation of $B$ for $l'\ge 4$.
If that is so, the latter separation is \emph{induced}
by the former one. The test uses the straightforward
subroutine PARTITION of the original method,
but nevertheless is rather complicated
since it exploits specific structural information concerning
$N^{i-1}$ and $N^i$. If an induced $(3|l')$-separation is
found, determine the corresponding 3-sum decomposition,
place the components into \BB, and go to Step~0.
In a rather complex process, retain additional
information that concerns the current sequence of nested
3-connected submatrices and, for special cases, facts about
induced decompositions. The retained information is used
later when the component
is removed from \BB\ and processed. It is precisely this
carryover of information that makes the relatively low
complexity of the original method possible.
If no induced $(3|l')$-separation is found, proceed as follows
depending on the classification of $N^i$ of Step~4: (1) If
$N^i$ is graphic, return to Step~2 to extend the current
sequence of nested 3-connected submatrices. (2) If $N^i$
is not graphic, equal to $B$, and equal to
one of $B^{10.1}$ or $B^{10.2}$, then
go to Step~0. (3) Otherwise, $N^i$ is not regular; declare
that $B$ of the original \BB\ is not regular, and stop.

\end{itemize}

%------------------------------------
\subsection{Implemented Method}
\label{sec:ImplementedMethod}
%------------------------------------

The main difference between the original and
the implemented method is replacement of the complicated Step~5
by a simpler, enumerating search that looks for induced
$(3|4)$-separations. The search is started
when (1) a 3-connected matrix $B$
not equal to $B^{10.1}$ or $B^{10.2}$ is at hand; (2) a sequence of
3-connected nested $N^1,\ldotspec$, $N^k = B$ has been found
where the length of $N^1$ satisfies $8\le s(N^1)\le 10$, and where
for $2\le i\le k$, the lengths of $N^{i}$ and $N^{i-1}$
satisfy $s(N^{i}) \le s(N^{i-1})+3$; and (3)
for some $1 \le j < k$ and all $1\le i\le j$, the matrices 
$N^i$ are graphic or cographic, while 
$N^{j+1}$ is not graphic and not cographic. 
Let $E^i$ be the union of the row and column index sets of 
the matrices $N^i$ of the sequence.

Since $N^{j+1}$ is not graphic and not cographic, $B$ is not
graphic and not cographic. If $B$ is regular, then
by Theorem~\ref{thm:Seymour}(v) there is a 
$(3|6)$-separation of $B$ that can be reduced to a
$(3|6)$-separation of $N^{j+1}$ by suitable reduction of
the index sets of the separation. Since 
$s(N^{i}) \le s(N^{i-1})+3$, the matrix $N^j$ has a 3-separation that
induces the $(3|6)$-separation of $B$. 
In the search described next,
we look instead for a less demanding $(3|4)$-separation of $B$
induced by a possibly deficient 3-separation of one of the
matrices $N^1,\ldotspec$, $N^j$.

We first generate all pairs $(T, E^1 \setminus T)$ where
(1) $T$ satisfies $|T| \leq |E^1 \setminus T|$, and
(2) $X_1 = X\cap T$, $X_2 = X\setminus X_1$,
$Y_1 = Y\cap T$, and $Y_2 = Y\setminus Y_1$ 
define a possibly deficient 3-separation of $N^1$.
The number of such pairs is bounded
by a constant since $|E^1| \leq 10$.
For $i = 2,\ldotspec$, $j$, we then
generate all pairs $(T,T\setminus E^i)$ where
(1) $T$ contains at least
one element of $E^i\setminus E^{i-1}$
and at most one element of $E^{i-1}$, and
(2) the pair defines a possibly deficient 3-separation
of $N^i$ analogously to the $N^1$ case.
Thus, the number of pairs for case $i$ is linear
in $|E^i|$, and overall a total of $O((m+n)^2)$ pairs
are produced.

The construction rules of the pairs obtained from the
sequence $N^1,\ldotspec$, $N^j$ assure validity
of the following claim. If $B$ is regular, then,
for some $i \le j$,
the possibly deficient 3-separation of $N^i$ corresponding to one
of the derived pairs induces a $(3|4)$-separation of $B$.

We use PARTITION of the original method to derive such an induced
$(3|4)$-separation and thus a 3-sum decomposition of $B$,
or to conclude that no such decomposition is possible.
PARTITION specifies that the smaller set $T$ of the input
pair satisfies $|T|\ge 3$, but with a trivial modification
the algorithm works just as well for the case $|T|=2$
arising here from deficient 3-separations.

Suppose a 3-sum decomposition is carried out. In the original
method, work done prior to that decomposition is used
when component matrices are processed. This approach
necessitates that tests for graphicness are accompanied
with a search for decompositions. The implemented
method avoids the complexity of that approach and simply
proceeds recursively after each decomposition. In the
simplified process, a sequence of nested 3-connected
submatrices $N^1,\ldotspec$, $N^k = B$ is determined,
and a test of
graphicness and cographicness is carried out for each
matrix of the sequence until either $B$ is determined to
be graphic or cographic, or a nongraphic and noncographic
$N^{j+1}$ is found. In the latter case, the above-described search
for a 3-sum decomposition either
produces such a decomposition or results in the conclusion
that the matrix is not regular.

We turn to the problem of signing the original $B$.
The original method does such signing
and then attempts to scale the signed version so that
it becomes the input matrix $A$.
In the implemented method, the signing of $A$ is taken into account
when $B$ is signed, and thus the scaling step is not needed.
Also, the signing is done at the very beginning since it
is quickly done and may already determine the input matrix $A$
to be non-t.u.

The overall run time of the implemented method is
$O((m+n)^5)$, since
(1) PARTITION has at most quadratic run time,
(2) a given nested sequence of 3-connected matrices
produces at most a quadratic number of input pairs for
PARTITION,
(3) there are at most a linear number of 3-sum
decompositions, and
(4) the remaining steps are easily done in $O((m+n)^5)$ time.

The user might want to obtain a certificate together with the answer.
A positive certificate consists of a tree whose inner nodes
correspond to 1-, 2-, and 3-sum decompositions, and whose
end nodes correspond to graphs and copies of $B^{10.1}$ and
$B^{10.1}$.

A negative certificate consists of a square matrix
whose determinant has absolute value of at least 2. More interesting
is a minimal violator that is not t.u., but all of whose proper
submatrices are t.u.
A simple strategy for finding a minimum violator recursively
removes a single row or column, tests for total unimodularity,
and adds the row or column back in if the submatrix turns out
to be t.u. Three ideas reduce the number of total unimodularity
tests during the search.

First, if the signing process determines $A$ to be non-t.u.
while $B$ was found to be regular, then the signing directly
determines a minimal violator having exactly two nonzeros in each
row and column; see Section~\ref{sec:DetailsofImplementation}.

Second, suppose that the method has stopped since the currently
processed matrix is nonregular.
Due to the structure of the decompositions,
that nonregular matrix is obtainable from the
original matrix $B$ by a sequence of pivots followed by
deletion of some rows and columns. The sequence of pivots and the
deletions are readily determined. If a deleted row or column was
never involved in a pivot, then its deletion from $B$
produces a smaller nonregular matrix. Hence, we carry out all such
deletions.

The third idea is based on the fact that the given matrix
may contain a number of minimal violators. This need not be so,
as shown by the construction of \cite{Truemper92} where the
matrices having exactly one minimal violator are produced.
Guessing that a number of minimal violators are present,
we remove 80\% of the rows or columns
of the nonregular $B$. If the resulting submatrix is nonregular,
recursion is used. Otherwise, we go back to the original matrix
and remove 40\% of the rows or columns. We go on by always
halving the amount until we reach a fixed threshold.
If the heuristic fails to produce a nonregular submatrix,
we apply the na\"{\i}ve strategy and
remove only a single row or column at a time.

We use a straightforward
matrix implementation with $O(1)$ indexed
access. For conceptual simplicity,
subroutines sometimes move certain submatrices of a given
matrix to the top left corner or carry out column/row
permutations or transposition. 
For efficient handling of these cases,
we have implemented various \emph{matrix proxies}.
The most important one is an object
that refers to another matrix,
but applies row and column permutations beforehand.
With the help of generic programming,
 it can be used as a usual matrix,
 but only needs $O(1)$ time to swap two rows or columns.

Below, the implemented method is called
\emph{Decomposition Test}, for short DT.
Method DT followed by construction of a minimal violator,
if applicable, is called DT\&V.

%-----------------------------------------------
\section{Enumerative Methods}
\label{sec:EnumerativeMethods}
%-----------------------------------------------

We also have implemented  two na\"{\i}ve enumerative tests
for comparison purposes. The first one tests the
square submatrices using the criterion of Camion 
\cite{Camion63}, according to which a matrix $A$ is t.u.
if and only if, for every square submatrix $A'$ of $A$
with even row and column sums, the sum of the entries
of $A'$ is divisible by 4.
{\em Submatrix Test}, for short ST, applies that test
to the square submatrices of $A$ in increasing order.

The second one is based on the characterization of
total unimodularity by Ghouila-Houri \cite{GhouilaHouri62},
which says that a matrix $A$ is t.u. if and only if
for each column submatrix $A'$ of $A$, there
is a $\{\pm1\}$ vector $x$ such that $A'x$ is a $\{0,\pm1\}$
vector. Algorithm {\em Column Enumeration}, for short CE,
carries out that test in straightforward fashion.

%-----------------------------------------------
\section{Computational Results}
\label{sec:ComputationalResults}
%-----------------------------------------------

We have implemented DT, DT\&V, ST, and CE in C\verb!++! and have
applied them to three matrix classes
using an AMD Opteron with 2.3\,GHz. In the tables below,
the run time is measured in seconds and 
omitted if less than $0.1$ sec.

The first class consists of randomly generated matrices. 
If the generation is
carried out without some care, then the resulting matrices
most likely contain a non-t.u. $2\times 2$ submatrix and typically
are detected to  be non-t.u. by the signing process. 
Thus, the methods ST, CE, and
DT would settle the cases very quickly. 
To prevent that trivial outcome, we generate
$\{0,1\}$ matrices $B$ randomly and 
apply the signing procedure to obtain 
$\{0,\pm1\}$ matrices $A$ which are then tested.
The net effect of this change is that
(1) ST and CE are less likely to terminate 
due to a $2\times 2$ non-t.u. submatrix, 
and (2)  DT and DT\&V never detect non-total unimodularity 
in the signing process and thus
always carry out the generally difficult regularity test of $B$.

Details of the generation of  the random matrices 
are as follows. For each $p = 2/3$, 1/2, 1/4, and 1/8,
and each $n = 200$, 400, and 800,  
we randomly select ten $n\times n$ $\{0,1\}$ matrices $B$ where
$p$ is the probability that a given entry receives the value 1.  
Then we apply the signing procedure to each $B$
to obtain $\{0,\pm1\}$ matrices $A$.
Thus, $A$ is t.u. if and only if $B$ is regular. 
To each such matrix $A$, the methods ST, CE, 
DT, and DT\&V are applied.
In Table~\ref{table:RandomMatrices} below, the results listed
for each pair of $p$ and $n$ are the geometric means 
of the run times for the ten matrices.

\begin{table}[h]
  \caption{Running times for random matrices}
  \centering\begin{tabular}{|c|c|r|r|r|r|}
    \hline
    Size & $p$ &
    \multicolumn{1}{c|}{ST} &
    \multicolumn{1}{c|}{CE} &
    \multicolumn{1}{c|}{DT} &
    \multicolumn{1}{c|}{DT\&V} \\
    \hline
    $200 \times 200$ & 2/3 & 185.4   &     &     & 0.1 \\
    $400 \times 400$ & 2/3 & 2948.8  & 0.1 & 0.3 & 0.4 \\
    $800 \times 800$ & 2/3 & $>3600$ & 0.2 & 0.7 & 0.9 \\
    \hline
    $200 \times 200$ & 1/2 & 28.2    &     &     & \\
    $400 \times 400$ & 1/2 & 218.3   & 0.1 & 0.1 & 0.1 \\
    $800 \times 800$ & 1/2 & $>3600$ & 0.3 & 0.7 & 0.7 \\
    \hline
    $200 \times 200$ & 1/4 & 0.1     &     & 0.1 & 0.1 \\
    $400 \times 400$ & 1/4 & 0.2     & 0.1 & 0.3 & 0.3 \\
    $800 \times 800$ & 1/4 & 1.0     & 0.3 & 0.4 & 0.5 \\
    \hline
    $200 \times 200$ & 1/8 & 0.1     &     & 0.1 & 0.1 \\
    $400 \times 400$ & 1/8 & 0.2     & 0.1 & 0.2 & 0.3 \\
    $800 \times 800$ & 1/8 & 0.5     & 0.2 & 0.6 & 2.0 \\
    \hline
  \end{tabular}
  \label{table:RandomMatrices}
\end{table}

It turns out that all matrices are non-t.u., which is no surprise. 
The impressive performance of method CE is due to the fact that 
non-total unimodularity can be proved using few columns.
Method DT also handles all cases well since it typically finds
small 3-connected nongraphic and noncographic submatrices that result
in few candidate pairs for induced (3|4)-separations.
Method ST is third in performance
and works well except for the cases with $p=2/3$ and
$n=400, 800$.
The methods DT and DT\&V have similar run times. This is due to the
fact that the heuristic for finding minimal violators 
described in Section~\ref{sec:ImplementedMethod} is very effective
for the matrix class constructed here. 
%Note that there are a few seemingly inconsistent cases where
%the run time listed for DT is a bit larger than that for DT\&V,
%which in reality is impossible 
%since DT\&V initially carries out DT. For example,
%for the $400\times 400$ case with $p=1/8$,
%DT is claimed to require 237.4 sec and DT\&V 237.0 sec. The minor
%inconsistencies are due to inaccuracies of the timing routine.
%Similar inconsistencies surface later 
%in Table~\ref{table:NetworkMatrices}, where the listed times
%vary slightly for DT and DT\&V, while actual times must be identical
%for the two methods.

The second set of test matrices is generated from randomly 
generated directed networks. According to Theorem~\ref{thm:Seymour},
these matrices and their transposes are in some sense
the main building blocks of t.u. matrices.
We generate the networks by
constructing Erd\"{o}s-R\'{e}nyi graphs 
$G(n,p)$ (see \cite{Gilbert59}) 
and compute matrices representing the corresponding graphic matroids.
The parameters $n$ and $p$ are chosen such that
the resulting matrix is connected and of suitable size.
For each matrix size listed below, we generate one instance.
These matrices are t.u. and therefore should be difficult 
for the enumerative
methods ST and CE.  At the same time, DT should perform
well since the method never needs to find a 3-sum
decomposition or locate a minimal violator.
The data of Table~\ref{table:NetworkMatrices} 
support these predictions.

\begin{table}[h]
 \caption{Running times for network matrices}
 \centering\begin{tabular}{|c|r|r|r|r|}
   \hline
   Size                & \multicolumn{1}{c|}{ST} &
\multicolumn{1}{c|}{CE} & \multicolumn{1}{c|}{DT} \\
   \hline
   $10 \times 10$      & 0.3     &         &     \\
   $12 \times 12$      & 1.8     &         &     \\
   $14 \times 14$      & 21.6    & 0.1     &     \\
   $16 \times 16$      & 311.8   & 1.0     &     \\
   $18 \times 18$      & $>3600$ & 8.3     &     \\
   $20 \times 20$      & $>3600$ & 83.1    &     \\
   $22 \times 22$      & $>3600$ & 709.6   &     \\
   $24 \times 24$      & $>3600$ & $>3600$ &     \\
   $50 \times 50$      & $>3600$ & $>3600$ &     \\
   $100 \times 100$    & $>3600$ & $>3600$ &     \\
   $200 \times 200$    & $>3600$ & $>3600$ & 0.1 \\
   $400 \times 400$    & $>3600$ & $>3600$ & 0.5 \\
   $800 \times 800$    & $>3600$ & $>3600$ & 5.3 \\
   \hline
 \end{tabular}
 \label{table:NetworkMatrices}
\end{table}

\noindent
The third class consists of matrices that contain
exactly one square submatrix whose determinant is
not equal to 0 or $\pm1$ and thus must be a
minimal violator. Each matrix is
constructed from a square matrix of odd order $n \ge 5$
where each row and each column has exactly two
1s, arranged in cycle fashion. In that matrix, two
adjacent rows are selected, and for any column
where those two rows have 0s, both 0s are replaced
by 1s. It is readily shown that any such matrix is
a minimal violator. Using results of
\cite{Truemper92}, the matrix is transformed by pivots
to one having exactly one violator of order 
$(n-1)/2$. Thus, as $n$ grows, effort of ST and CE must
grow exponentially, and even for modest values of $n$ 
the two methods should
be unable to decide total unimodularity. On the other
hand, DT should be able to process these
difficult cases with reasonable efficiency. Finally,
DT\&V should require significant additional effort
beyond that for DT since each matrix has just one
minimal violator, and that submatrix is relatively large.
Table~\ref{table:MinimalViolators}
confirms these predictions.

\begin{table}[h]
 \caption{Running times for odd-cycle matrices}
 \centering\begin{tabular}{|c|r|r|r|r|}
   \hline
   Size                & \multicolumn{1}{c|}{ST} &
\multicolumn{1}{c|}{CE} & \multicolumn{1}{c|}{DT} &
\multicolumn{1}{c|}{DT\&V} \\
   \hline
   $11 \times 11$      & 0.7     &         &         & \\
   $13 \times 13$      & 5.7     &         &         & \\
   $15 \times 15$      & 79.3    & 0.6     &         & \\
   $17 \times 17$      & 1201.8  & 6.0     &         & 0.1 \\
   $19 \times 19$      & $>3600$ & 59.1    &         & 0.1 \\
   $21 \times 21$      & $>3600$ & 569.2   &         & 0.2 \\
   $23 \times 23$      & $>3600$ & $>3600$ &         & 0.3 \\
   $51 \times 51$      & $>3600$ & $>3600$ & 0.6     & 3.2 \\
   $101 \times 101$    & $>3600$ & $>3600$ & 3.4     & 170.7 \\
   $151 \times 151$    & $>3600$ & $>3600$ & 13.1    & 1302.5 \\
   $201 \times 201$    & $>3600$ & $>3600$ & 32.1    & 3490.1 \\
   $301 \times 301$    & $>3600$ & $>3600$ & 137.2   & $>3600$ \\
   $401 \times 401$    & $>3600$ & $>3600$ & 476.3   & $>3600$ \\
   $501 \times 501$    & $>3600$ & $>3600$ & 1457.1  & $>3600$ \\
   $601 \times 601$    & $>3600$ & $>3600$ & 2422.0  & $>3600$ \\
   $701 \times 701$    & $>3600$ & $>3600$ & $>3600$ & $>3600$ \\
   \hline
 \end{tabular}
 \label{table:MinimalViolators}
\end{table}

%-----------------------------------------------
\section{Details of Implementation }
\label{sec:DetailsofImplementation}
%-----------------------------------------------

This section provides implementation details.
We skip the trivial enumerative tests
ST and CE, and also omit steps
that are identical to their counterparts 
in the original method. As an easy means of differentiation, we use
the names of the procedures of the original method in
capital letters, just as done in \cite{Truemper90}, and
employ lower-case names for the subroutines introduced here.
Bold fonts for both types of names help set them apart from
the text.

The discussion proceeds in a top-down manner. Thus, we start
with the main routine, which tests for total unimodularity.

\algorithmName{{\bf is\_totally\_unimodular}}
\algorithmInput{${m \times n}$ matrix $A$ with $\{0, \pm 1\}$ entries.}
\algorithmOutput{\True{} if $A$ is t.u., and \False\ otherwise.}
\algorithmComplexity{$O((m + n)^5)$ time and $O(m \cdot n)$ space.}
\algorithmBegin
\item
  Call {\bf sign\_matrix} with $A$ as input. If the subroutine
declares that $A$ has been modified, return \False.
\item
  Call {\bf decompose\_matrix} with the binary version $B$ of $A$
  and an empty sequence of nested 3-connected submatrices
  as input. Return the \True/\False\ output of the subroutine.
\algorithmEnd
\algorithmLine

The signing procedure described next recursively modifies the entries 
of the $\{0,\pm1\}$ input matrix $A$ so that certain 
minimal submatrices $V$
with exactly two nonzeros in each row and column
and with entries summing to $2 \pmod 4$ have the entries sum to
$0 \pmod 4$ after the signing. In the context of the main routine
{\bf is\_totally\_unimodular}, the matrix $A$ is t.u. if and only if
the signing routine does not change any entry of $A$ and
the binary version $B$ of $A$ is regular. 
On the other hand, if the 
signing routine does change at least one such entry, then
the subroutine can
be stopped when the first such change is to be made. The matrix
$V$ on hand at that time is a minimal non-t.u matrix.

Recall that the signing procedure is used in the construction
of the first class of test matrices, where each randomly generated
matrix $B$ is signed by the signing procedure to obtain a matrix
$A$. 
When methods DT and DT\&V later process each such test matrix $A$,
the signing procedure does not alter $A$, 
and thus regularity of $B$ is tested. The outcome
of the latter test then decides whether $A$ is t.u.

\algorithmName{{\bf sign\_matrix}}\label{alg:sign_matrix}
\algorithmInput{${m \times n}$ matrix $A$ with $\{0, \pm 1\}$ entries.}
\algorithmOutput{Either: ``$A$ has been modified.'' 
                 Or: ``$A$ is unchanged.''}
\algorithmComplexity{$O( m \cdot n^2)$ time and $O(m \cdot n)$ space.}
\algorithmBegin
\item\label{algo:DT:sign:init}
  If $A$ is a zero matrix, output ``$A$ is unchanged'', and stop.
  Otherwise, let $y$ be the index of an arbitrarily selected
  nonzero column of $A$. Define $X$ to be the index set of
  the rows of $A$ containing the nonzeros of the selected column,
  and initialize $Y = \{y\}$.
\item\label{algo:DT:sign:test}
  If $Y$ is equal to the column index 
  set of $A$, output ``$A$ has been modified'' if
  during any iteration in Step~\ref{algo:DT:sign:BFS} 
  an entry was changed, and output
  ``$A$ is unchanged'' otherwise; stop.
\item
  If there is a column index $y\notin Y$ for which the column vector
  $A_{X,y}$ is nonzero, select one such $y$,
  and go to Step~\ref{algo:DT:sign:BFS}. Otherwise, select
  any $y\notin Y$, and go to Step~\ref{algo:DT:sign:expand}.
\item \label{algo:DT:sign:BFS}
  Let $A_{x_j,y}$, $j=0,\ldotspec$, $s$, be the nonzero entries
  of the column vector $A_{X,y}$. If the vector has just one
  nonzero entry and 
  thus $s=0$, go to Step~\ref{algo:DT:sign:expand}.
  Otherwise, do a breadth-first-search on $\BG\left(A_{X,Y}\right)$
  to find shortest paths from $x_0$ to every other $x_j$, 
  $j=1,\ldotspec$, $s$.

  \quad For $i = 1,\ldotspec$, $s$, do the following steps. Going from
  $x_i$ to $x_0$ along the given shortest path, let the first
  $x_j$, $j\neq i$ encountered have index $j = p_i$.
  The path segment from $x_{p_i}$ to $x_i$
  together with the edges $(x_i, y)$ and $(x_{p_i}, y)$
  form a chordless cycle. In $A$, that cycle corresponds to a
  square submatrix $V$ with exactly two nonzero entries in each row
  and column. If the entries of that submatrix sum to
  $2 \pmod 4$, flip the sign of $A_{x_j,y}$.
\item \label{algo:DT:sign:expand}
  Add $y$ to $Y$. 
  Add the row indices $x$ for which $A_{x,y} \neq 0$ to $X$.
  Go to Step~\ref{algo:DT:sign:test}. 
\algorithmEnd
\algorithmRemark{Every iteration for given index $y$ 
is done in $O(m \cdot n)$ time.}
\algorithmLine

The next subroutine is the top procedure for the regularity test
of $B$. Given $B$, either (1) the subroutine decomposes $B$ in
a 1- or 2-sum decomposition and then invokes recursion 
for the components; or (2) it determines $B$ to be isomorphic
to $B^{10.1}$ or $B^{10.2}$ of Fig.~\ref{fig:MatricesforR10}
and thus to be regular;
 or (3) it constructs a sequence of nested
3-connected submatrices $N^1,\ldotspec$, $N^k$ which are tested
for graphicness and cographicness; if all submatrices turn
out to be graphic or cographic, $B$ is regular; or (4) it uses
the sequence $N^1,\ldotspec$, $N^k$ to decompose
$B$ in a 3-sum decomposition and then invokes recursion for
the components of the decomposition;  or, if none of (1)-(4)
apply, (5) it concludes that $B$ is not regular.

\algorithmName{{\bf decompose\_matrix}}
\algorithmInput{Binary ${m \times n}$ matrix $B$  and
possibly empty sequence of
nested 3-connected submatrices $N^1,\ldotspec$, $N^k$.}
\algorithmOutput{\True{} if $B$ is regular, and \False{} otherwise.}
\algorithmComplexity{$O( (m+n)^5)$ time and $O(m \cdot n)$ space.}
\algorithmBegin
\item
  If $m < 3$ or $n < 3$, then $B$ is regular; return \True{}, and stop.
\item
  If the  input sequence of nested 3-connected submatrices is 
  nonempty, go to Step~\ref{algo:DT:decompose:extend}.
  Otherwise, go to Step~\ref{algo:DT:decompose:W3}.
\item \label{algo:DT:decompose:W3}
  Call {\bf FIND-B}. 
  If it finds a 1- or 2-separation,
  decompose $B$ into $B_1$ and $B_2$
  according to the separation, and call
  {\bf decompose\_matrix} with each of them
  and the empty sequence as input.
  Return \True{} if both
  calls return \True{}, return \False{} otherwise, and stop.
  
  \quad If {\bf FIND-B} does not find a 1- or 2-separation, it has
  identified a submatrix $N^1$ that represents the graph $W_3$. 
  Initialize $k = 1$, and go to Step~\ref{algo:DT:decompose:extend}.
\item \label{algo:DT:decompose:extend}
  If $N^k = B$, go to Step~\ref{algo:DT:check:B10}.
  Otherwise, call {\bf EXTEND-B}, and try 
  to find a submatrix $N^{k+1}$ of $B$ which contains $N^k$.

  \quad If {\bf EXTEND-B} succeeds, increment $k$ by 1, and repeat 
  this step. Otherwise, the subroutine returns a
  1- or 2-separation of $B$ into $B_1$ and $B_2$ where
  $B_1$ contains a 3-connected submatrix sequence that is isomorphic
  to $N^1,\ldotspec$, $N^k$. 
  Relabel that sequence as $N^1,\ldotspec$, $N^k$.
  Call {\bf decompose\_matrix} twice: once with
  $B_1$ and $N^1,\ldotspec$, $N^k$, and the second time
  with $B_2$ and the empty sequence.
  Return \True{} if both calls return \True, return \FALSE{}
  otherwise, and stop.
\item \label{algo:DT:check:B10}
  If $B$ is $5 \times 5$, test whether
  $\BG(B)$ is isomorphic to $\BG(B^{10.1})$ 
  or $\BG(B^{10.2})$, where $B^{10.1}$ and $B^{10.2}$
  are shown in Fig.~\ref{fig:MatricesforR10}.
  The isomorphism test is trivial due to the size and special
  structure of $B^{10.1}$ and $B^{10.2}$.
  In the affirmative case, return \True, and stop.
\item \label{algo:DT:decompose:graph}
  Call {\bf test\_graphicness} to decide graphicness of $B$
  via the sequence $N^1,\ldotspec$, $N^k$. If the output is 
  ``$B$ is graphic,'' return \True{}, and stop. Otherwise,
  call {\bf test\_graphicness} once more to decide graphicness
  of $\transpose{B}$ 
  via $\transpose{N^1},\ldotspec$, $\transpose{N^k}$.
  If the output is ``$B^t$ is graphic,'' return \True, and stop.
\item
  In Step \ref{algo:DT:decompose:graph}, each of the two calls of
  {\bf test\_graphicness} returned an index. Let $j$
  the larger of the two returned indices. Thus,
  $N^1,\ldotspec$, $N^j$ are graphic or cographic, while
  $N^{j+1}$ is not graphic and not cographic.

  \quad For every $i = 1,\ldotspec$, $j$,
  define pairs $(T, E^i \setminus T)$ corresponding to
  possibly deficient 3-separations of $N^i$, as stated in 
  Section~\ref{sec:ImplementedMethod}. 
  For each such pair, call {\bf PARTITION} 
  to test whether the pair induces
  a $(3|4)$-separation of $B$. As soon as
  an induced decomposition of $B$ is detected,
  say into $B_1$ and $B_2$,
  call {\bf decompose\_matrix} twice, once with
  $B^1$ and the empty sequence, and the second time with
  $B^2$ and the empty sequence. Return
  \True{} if both calls return \True{}, return \FALSE{}
  otherwise, and stop.

  \quad If none of the pairs defined 
  for the matrices $N^1,\ldotspec$, $N^j$
  induces a $(3|4)$-separation of $B$, return \False, and stop.
\algorithmEnd
\algorithmRemark{The implementations of {\bf FIND-B} and
  {\bf EXTEND-B} run in $O((m+n)^3)$ time,
  while that of {\bf PARTITION} runs in $O(m \cdot n)$ time.}
\algorithmLine

Subroutine {\bf test\_graphicness} is described next.
It calls {\bf TEST-C} repeatedly to extend the
graph of a submatrix $N^i$ to the graph of the next
larger submatrix $N^{i+1}$.
\textbf{TEST-C} contains a minor error arising from two
special cases. The code uses the amended version.

\algorithmName{{\bf test\_graphicness}}
\algorithmInput{Binary $m\times n$ matrix $B$
and nonempty sequence $N^1,\ldotspec$, $N^k = B$ 
of nested 3-connected submatrices.}
\algorithmOutput{Either: ``$B$ is graphic'' together with
 the corresponding graph $G$.
 Or: ``$B$ is not graphic'' and the index of the largest graphic
 $N^i$.}
\algorithmComplexity{$O( (m+n)^3)$ time and $O(m \cdot n)$ space.}
\algorithmBegin
\item
  Define $G^1$ to be the wheel graph $W_3$, which represents $N^1$.
\item
  For $i=1,\ldotspec$, $k-1$, call {\bf TEST-C} to attempt extension
  of the graph $G^i$ for $N^i$ to a graph
  $G^{i+1}$ for $N^{i+1}$. If such an extension is not possible
  for some $N^i$, return ``$B$ is not graphic'' and index $i$,
  and stop.
  Otherwise, return ``$B$ is graphic'' together with the graph $G^k$
  for $N^k$, and stop.
\algorithmEnd
\algorithmRemark{The implementation of {\bf TEST-C} 
  runs in $O(m \cdot n)$.}
\algorithmLine

If $A$ is not t.u., a minimal violator is found as described in
Section~\ref{sec:ImplementedMethod}. 
Finally, tests for unimodularity and strong unimodularity have
been implemented in straightforward fashion 
using Theorem~\ref{thm:unimodularity}.
In the unimodularity test of an integer matrix $A$,
Gaussian elimination selects a basis $C$ of $A$ and
solves the equation $CX = A$ for $X$. Computation of the
Smith Normal Form \cite{Smith1861} settles whether
$\gcd_i\det C^i = 1$, and {\bf is\_totally\_unimodular}
decides if $X$ is t.u. The matrix $A$ is unimodular if and
only if $\gcd_i\det C^i = 1$ and $X$ is t.u.
Strong unimodularity of $A$ is decided by testing if
both $A$ and $A^t$ are unimodular.

%-----------------------------------------------
\section{Summary}
\label{sec:Summary}
%----------------------------------------------

The work reported here was motivated by the practical need
for an effective computer program testing
total unimodularity, strong unimodularity, and unimodularity.
We decided to use the algorithm of \cite{Truemper90} for the total
unimodularity test, which for $m\times n$ matrices 
has $O((m+n)^3)$ complexity.
Exact implementation of that algorithm is rather daunting,
so we used a somewhat simplified version with
$O((m+n)^5)$ complexity. Computational tests 
indicate that matrices with several
hundred rows and columns should be handled in reasonable time.
For the processing of larger matrices, the implementation could
be refined using the parts of the algorithm of 
\cite{Truemper90} that have been simplified here. 
An intermediate remedy could also be implemented where
the time-consuming search for decompositions is handled
by parallel processors. The present program is purposely
structured to simplify such a shift to parallel computation.

%\bibliographystyle{abbrv}
\bibliography{references}

\begin{thebibliography}{10}

\bibitem{Bixby80}
R.~E. Bixby and W.~H. Cunningham.
\newblock Converting linear programs to network problems.
\newblock {\em Math. Op. Res.}, 5:321--357, 1980.

\bibitem{Bixby86}
R.~E. Bixby, W.~H. Cunningham, and R.~Rajan.
\newblock A decomposition algorithm for matroids.
\newblock Technical report, Rice University, 1986.

\bibitem{Boostlicense}
Boost Software License\\ \texttt{http://www.boost.org/LICENSE\_1\_0.txt}.

\bibitem{Camion63}
P.~Camion.
\newblock {\em Matrices Totalement Unimodulaire et Probl\`emes Combinatoires}.
\newblock PhD thesis, Universit\'e Libre de Bruxelles, Bruxelles, 1963.

\bibitem{Cunningham78}
W.~H. Cunningham and J.~Edmonds.
\newblock Decomposition of linear systems, (unpublished).
\newblock 1965.

\bibitem{Edmonds65}
J.~Edmonds.
\newblock Minimum partition of a matroid into independet subsets.
\newblock {\em J. Res. Nat. Bur. Std. (B)}, 69:67--72, 1965.

\bibitem{Fujishige80}
S.~Fujishige.
\newblock An efficient $pq$-graph algorithm for solving the graph-realization
  problem.
\newblock {\em J. Computer and Systems Sciences}, 21:63--86, 1980.

\bibitem{GhouilaHouri62}
A.~Ghouila-Houri.
\newblock Caracterisation des matrices totalement unimodulaires.
\newblock {\em C.R. Acad. Sci. Paris}, 254:1192--1194, 1962.

\bibitem{Gilbert59}
E.~N. Gilbert.
\newblock Random graphs.
\newblock {\em Annals of Mathematical Statistics}, 30:1141--1144, 1959.

\bibitem{Hoffman56}
A.~J. Hoffman and J.~B. Kruskal.
\newblock Integral boundary points of convex polyhedra.
\newblock {\em Annals Math. Studies}, 38:223--246, 1956.

\bibitem{Hoffman78}
A.~J. Hoffman and R.~Oppenheim.
\newblock Local unimodularity in the matching polytope.
\newblock {\em Annals of Discrete Mathematics}, 2:201--209, 1978.

\bibitem{Seymour80}
P.~D. Seymour.
\newblock Decomposition of regular matroids.
\newblock {\em J. Comb. Theory, Ser. B}, 28:305--359, 1980.

\bibitem{Smith1861}
H.~J.~S. Smith.
\newblock On systems of linear indeterminate equations and congruences.
\newblock {\em Philos. Trans. Roy. Soc. London}, 151:293--326, 1861-1862.

\bibitem{Software}
Software available at:\\
  \texttt{http://www.utdallas.edu/{\url{~}}klaus/TUtest/}\\
  \texttt{http://www.math.uni-magdeburg.de/{\url{~}}walter/TUtest/}.

\bibitem{Truemper78}
K.~Truemper.
\newblock Algebraic characterizations of unimodular matrices.
\newblock {\em SIAM J. Appl. Math.}, 35:328--332, 1978.

\bibitem{Truemper80}
K.~Truemper.
\newblock Complement total unimodularity.
\newblock {\em Lin. Alg. Appl.}, 30:77--92, 1980.

\bibitem{Truemper90}
K.~Truemper.
\newblock A decomposition theory for matroids. {V}. {T}esting of matrix total
  unimodularity.
\newblock {\em J. Comb. Theory, Ser. B}, 49:241--281, 1990.

\bibitem{Truemper92}
K.~Truemper.
\newblock A decomposition theory for matroids. {VII}. {A}nalysis of minimal
  violation matrices.
\newblock {\em J. Comb. Theory, Ser. B}, 55:302--335, 1992.

\bibitem{Truemper98}
K.~Truemper.
\newblock {\em Matroid Decomposition (Revised Edition)}.
\newblock Leibniz Company, 1998.

\bibitem{TruemperRC78}
K.~Truemper and R.~Chandrasekaran.
\newblock Local unimodularity of matrix-vector pairs.
\newblock {\em Lin. Alg. Appl.}, 22:65--78, 1978.

\bibitem{Tutte58}
W.~T. Tutte.
\newblock A homotopy theorem for matroids i, ii.
\newblock {\em Trans. Amer. Math. Soc.}, 88:527--552, 1958.

\bibitem{Whitney33}
H.~Whitney.
\newblock 2-isomorphic graphs.
\newblock {\em Am. J. Math.}, 55:245--254, 1933.

\end{thebibliography}
% haystack.bib is the name of the Bibliography in this case

%\input haystack.bbl
% For submission:
% 1. Manually add tilde to 'Acuna' reference in bbl file
%    Convert 'minsat' to 'MINSAT'
% 2. Copy into the position '\input haystack.bbl' the
%    bbl file. Comment out '\bibliography{haystack}'

\end{document}